\RequirePackage{ifpdf}
\ifpdf % We are running pdfTeX in pdf mode
\documentclass[pdftex]{sigma}
\else
\documentclass{sigma}
\fi

\numberwithin{equation}{section}

\newtheorem{Theorem}{Theorem}[section]
\newtheorem{Lemma}[Theorem]{Lemma}
\newtheorem{Proposition}[Theorem]{Proposition}
{\theoremstyle{definition}
\newtheorem{Definition}[Theorem]{Definition}
\newtheorem{Remark}[Theorem]{Remark}
\newtheorem{Note}[Theorem]{Note}

\begin{document}

\allowdisplaybreaks

\renewcommand{\PaperNumber}{071}

\FirstPageHeading

\ShortArticleName{Inversion of the Dual Dunkl--Sonine Transform on $\mathbb{R}$}

\ArticleName{Inversion of the Dual Dunkl--Sonine Transform on $\boldsymbol{\mathbb{R}}$\\ Using Dunkl Wavelets}

\Author{Mohamed Ali MOUROU}

\AuthorNameForHeading{M.A. Mourou}
\Address{Department of Mathematics, Faculty of Sciences, Al-Jouf University,\\
P.O.~Box 2014, Al-Jouf, Skaka, Saudi Arabia}
\Email{\href{mailto:mohamed_ali.mourou@yahoo.fr}{mohamed\_ali.mourou@yahoo.fr}}

\ArticleDates{Received March 02, 2009, in f\/inal form July 04, 2009;  Published online July 14, 2009}

\Abstract{We prove a Calder\'{o}n reproducing formula for the Dunkl continuous wavelet transform on $\mathbb{R}$.
We apply this result to derive new inversion formulas for the dual Dunkl--Sonine integral transform.}


\Classification{42B20; 42C15; 44A15; 44A35}

\section{Introduction}

The one-dimensional Dunkl kernel $e_\gamma$, $\gamma>-1/2$, is def\/ined by \[e_{\gamma}(z) = j_{\gamma}(iz) +
\frac{z}{2(\gamma + 1)} j_{\gamma + 1} (iz), \qquad  z \in \mathbb{C} ,\]
where
\[
j_\gamma(z) = \Gamma(\gamma + 1)
\sum^\infty_{n=0} \frac{(-1)^n\; (z/2)^{2n}}{n!\; \Gamma(n + \gamma+1)} %\quad (z \in \mathbb{C}),
\] is the
normalized spherical Bessel function of index $\gamma$.   It is well-known (see~\cite{Dunkl1991}) that the functions
$e_\gamma(\lambda \cdot)$, $\lambda \in \mathbb{C}$, are solutions of the dif\/ferential-dif\/ference equation
\[
\Lambda_\gamma u = \lambda u, \qquad u(0) = 1,
\]
where \[\Lambda_\gamma f(x) = f'(x) + \left(\gamma + \frac{1}{2}\right)
\frac{f(x) - f(-x)}{x}
\] is the  Dunkl operator with parameter $\gamma+1/2$ associated with the ref\/lection grour
$\mathbb{Z}_2$ on $\mathbb{R}$. Those operators were introduced and studied by Dunkl \cite{Dunkl1989, Dunkl1991, Dunkl1992} in connection with
a generalization of the classical theory of spherical harmonics. Besides its mathematical interest, the Dunkl
operator $\Lambda_{\alpha}$ has quantum-mechanical applications; it is naturally involved  in the study of
one-dimensional harmonic oscillators governed by Wigner's commutation rules \cite{Kamefuchi, Rosenblum, Yang}.

It is known, see for example \cite{Soltani2008, Xu}, that the Dunkl kernels on $\mathbb{R}$ possess the following Sonine type
integral representation \begin{gather}
e_{\beta}(\lambda
x)=\int_{-|x|}^{|x|}\mathcal{K}_{\alpha,\beta}(x,y)\,e_{\alpha}(\lambda y)\,|y|^{2\alpha+1}\,dy, \qquad
 \lambda\in\mathbb{C}, \qquad x\neq 0, \label{eq1}
 \end{gather}
 where
\begin{gather}\mathcal{K}_{\alpha,\beta}(x,y):= \left\{
\begin{array}{ll}
 a_{\alpha,\beta}\,\mbox{sgn}(x)\;(x+y)\,\displaystyle{\frac{\left(x^2-y^2\right)^{\beta-\alpha-1}}{|x|^{2\beta+1}}}
 &\mbox{ if } |y|<|x|,\\\\
 0&\mbox{ if } |y|\geq|x|,
 \end{array}\right.\label{eq2}
 \end{gather}
 with
 $\beta>\alpha>-1/2$, and \[ a_{\alpha,\beta}:=\frac{\Gamma(\beta+1)}{\Gamma(\alpha
+1)\,\Gamma(\beta-\alpha)}.
\]

Def\/ine the Dunkl--Sonine integral transform $\mathcal{X}_{\alpha,\beta}$ and its dual
$^t\!\mathcal{X}_{\alpha,\beta}$, respectively, by
\begin{gather*}
\mathcal{X}_{\alpha,\beta}f(x)=\int_{-|x|}^{|x|}\mathcal{K}_{\alpha,\beta}(x,y)\,f(y) \,|y|^{2\alpha+1}\,dy,\\
 ^t\!\mathcal{X}_{\alpha,\beta}f(y)=\int_{|x|\geq |y|}\mathcal{K}_{\alpha,\beta}(x,y)\,f(x)
\,|x|^{2\beta+1}\,dx.
\end{gather*}

Soltani has showed in \cite{Soltani2008} that the dual Dunkl--Sonine integral transform
 $^t\!\mathcal{X}_{\alpha,\beta}$ is a transmutation operator between $\Lambda_\alpha$ and $\Lambda_\beta$
  on the Schwartz space $\mathcal{S}(\mathbb{R})$, i.e., it is an automorphism of~$\mathcal{S}(\mathbb{R})$ satisfying
the intertwining relation
\[^t\!\mathcal{X}_{\alpha,\beta}\,\Lambda_{\beta}\,f=\Lambda_{\alpha}\,^t\!\mathcal{X}_{\alpha,\beta}\,f, \qquad
f\in \mathcal{S}(\mathbb{R}).\]The same author \cite{Soltani2008} has obtained inversion formulas for the transform
$^t\!\mathcal{X}_{\alpha,\beta}$ involving pseudo-dif\/ferential-dif\/ference operators and only valid on a restricted
subspace of $\mathcal{S}(\mathbb{R})$.

The purpose of this paper is to investigate the use of Dunkl wavelets (see \cite{Jouini}) to derive a new inversion of the
dual Dunkl--Sonine transform on  some Lebesgue spaces. For other applications of wavelet type transforms to inverse
problems we refer the reader to \cite{Mourou1998,Trimeche1998} and the references therein.

The content of this article is as follows. In Section~\ref{section2} we recall some basic harmonic analysis results related to
the Dunkl operator. In Section~\ref{section3} we list some basic properties of the Dunkl--Sonine integral trnsform and its dual.
In Section~\ref{section4} we give the def\/inition of the Dunkl continuous wavelet transform and we establish for this transform
a Calder\'{o}n formula. By combining the results of the two previous sections, we obtain in Section~\ref{section5} two new
inversion formulas for the dual Dunkl--Sonine
 integral transform.

\section{Preliminaries}\label{section2}

\begin{Note} \label{note1} Throughout this section assume $\gamma
> - 1/2$. Def\/ine $L^p(\mathbb{R},
|x|^{2\gamma+1} dx)$, $1 \leq p \leq \infty$, as the class of measurable  functions $f$ on $\mathbb{R}$ for which
$||f||_{p, \gamma} < \infty$, where \[
||f||_{p, \gamma} = \left(\int_{\mathbb{R}} |f(x)|^p |x|^{2\gamma+1}
dx\right)^{1/p},\qquad \mbox{if} \quad p < \infty,\] and $||f||_{\infty, \gamma} =||f||_{\infty}= {\rm ess\,sup}_{x
\in \mathbb{R}}|f(x)|$. $\mathcal{S}(\mathbb{R})$ stands for the usual Schwartz space.
\end{Note}

The Dunkl transform  of order $\gamma + 1/2$ on $\mathbb{R}$ is def\/ined for a function $f$ in $ L^1(\mathbb{R},
|x|^{2\gamma+1} dx)$ by
\begin{gather}{\cal F}_\gamma f(\lambda) = \int_{\mathbb{R}} f(x)\,e_\gamma(- i\lambda
x)\,|x|^{2\gamma+1}dx,\qquad \lambda \in \mathbb{R}. \label{eq3}
\end{gather}

\begin{Remark}\label{remark1} It is known that
 the Dunkl transform ${\cal F}_\gamma$ maps continuously and
injectively \linebreak $ L^1(\mathbb{R}, |x|^{2\gamma+1} dx)$ into the space $\mathcal{C}_0(\mathbb{R})$ (of continuous
functions on $\mathbb{R}$ vanishing at inf\/inity).
\end{Remark}

Two  standard results about the Dunkl transform ${\cal F}_\gamma$ are as follows.
\begin{Theorem}[see~\cite{De Jeu}]\label{theorem1}  \qquad

\begin{enumerate}\itemsep=0pt
\item[$(i)$] For every $f \in L^1 \cap L^2(\mathbb{R}, |x|^{2\gamma+1}
dx)$ we have the  Plancherel formula
\[\int_{\mathbb{R}} |f(x)|^2 |x|^{2\gamma+1}dx = m_\gamma\int_{\mathbb{R}}|{\cal F}_\gamma f(\lambda)|^2
|\lambda|^{2\gamma+1}d\lambda,
\]
where
\begin{gather}
m_\gamma=\frac{1}{2^{2\gamma+2}(\Gamma(\gamma+1))^2}.\label{eq4}
\end{gather}

\item[$(ii)$] The Dunkl transform ${\cal F}_\alpha$ extends uniquely to an isometric isomorphism from \linebreak $L^2(\mathbb{R},
|x|^{2\gamma+1}dx)$ onto $L^2(\mathbb{R}, m_\gamma |\lambda|^{2\gamma+1}d\lambda)$. The inverse transform is given
by \[{\cal F}^{-1}_\gamma g(x) = m_\gamma\int_{\mathbb{R}} g(\lambda) e_\gamma(i\lambda x)
|\lambda|^{2\gamma+1}d\lambda,
\] where the integral converges in $L^2(\mathbb{R}, |x|^{2\gamma+1} dx)$.
\end{enumerate}
\end{Theorem}

\begin{Theorem}[see~\cite{De Jeu}]\label{theorem2}  The Dunkl transform ${\cal F}_\alpha$ is an automorphism of
$\mathcal{S}(\mathbb{R})$.
\end{Theorem}

An outstanding result about Dunkl kernels on $\mathbb{R}$ (see~\cite{Rosler}) is the product formula
\[
e_\gamma(\lambda x)
e_\gamma(\lambda y) =T^x_\gamma \left(e_\gamma(\lambda \cdot)\right)(y),\qquad \lambda\in\mathbb{C}, \qquad x,y\in
\mathbb{R},
\] where $T^x_\gamma$ stand for  the  Dunkl translation operators def\/ined by
\begin{gather}
T^x_\gamma f(y) =
\frac{1}{2}\int_{-1}^1
f\left(\sqrt{x^2+y^2-2xyt}\right)\left(1+\frac{x-y}{\sqrt{x^2+y^2-2xyt}}\right)W_\gamma(t)dt\nonumber\\
\phantom{T^x_\gamma f(y) =}{} +\frac{1}{2}\int_{-1}^1
f\left(-\sqrt{x^2+y^2-2xyt}\right)\left(1-\frac{x-y}{\sqrt{x^2+y^2-2xyt}}\right)W_\gamma(t)dt,\label{eq5}
\end{gather}
with \[
W_\gamma(t)=\frac{\Gamma(\gamma+1)}{\sqrt{\pi}\,\Gamma(\gamma+1/2)}(1+t)\left(1-t^2\right)^{\gamma-1/2}.
\]

The Dunkl convolution of two functions $f$, $g$ on $\mathbb{R}$ is def\/ined by the relation
 \begin{gather}
 f
\ast_\gamma g(x) = \int_{\mathbb{R}} T^x_\gamma f(-y) g(y)|y|^{2\gamma+1}dy.\label{eq6}
\end{gather}
\begin{Proposition}[see \cite{Soltani2004}] \label{proposition1}\qquad

\begin{enumerate}\itemsep=0pt
\item[$(i)$] Let $p, q, r\in \![1,\infty]$ such that $\frac{1}{p}+\frac{1}{q}-1=\frac{1}{r}$.
If $f \in\! L^p(\mathbb{R}, |x|^{2\gamma+1} dx)$ and $g \in\! L^q(\mathbb{R}, |x|^{2\gamma+1} dx)$, then $f
\ast_\gamma g \in L^r (\mathbb{R}, |x|^{2\gamma+1} dx)$ and
\begin{gather}||f \ast_\gamma g||_{r, \gamma} \leq
4  ||f||_{p, \gamma}  ||g||_{q, \gamma}.\label{eq7}\end{gather}

\item[$(ii)$] For $f \in L^1 (\mathbb{R}, |x|^{2\gamma + 1}dx)$ and $g \in L^p (\mathbb{R}, |x|^{2\gamma + 1}dx)$, $p = 1$ or
$2$, we have
\begin{gather}
{\cal F}_\gamma (f \ast_\gamma g) = {\cal F}_\gamma f  {\cal F}_\gamma g. \label{eq8}
\end{gather}
\end{enumerate}
\end{Proposition}

For more details about harmonic analysis related to the Dunkl operator on $\mathbb{R}$ the reader is referred, for
example, to \cite{Mourou2002,Mourou2003}.

\section[The dual Dunkl-Sonine  integral transform]{The dual Dunkl--Sonine  integral transform}\label{section3}

Throughout this section assume $\beta>\alpha>-1/2$.

\begin{Definition}[see \cite{Soltani2008}]\label{definition1}
The dual Dunkl--Sonine
integral transform $^t\!\mathcal{X}_{\alpha,\beta}$ is def\/ined for smooth functions on  $\mathbb{R}$ by
\begin{gather}^t\!\mathcal{X}_{\alpha,\beta}f(y):=\int_{|x|\geq |y|}\mathcal{K}_{\alpha,\beta}(x,y) f(x)
 |x|^{2\beta+1}\,dx,\qquad y\in \mathbb{R} ,\label{eq9}
\end{gather}
where $\mathcal{K}_{\alpha,\beta}$ is the
kernel given by~\eqref{eq2}. \end{Definition}

 \begin{Remark}\label{remark2} Clearly, if ${\rm supp}\, (f) \subset [-a,a]$  then ${\rm supp} \left(
^t\!\mathcal{X}_{\alpha,\beta}f\right)\subset [-a,a]$.
\end{Remark}

The next statement provides formulas relating harmonic analysis tools tied to $\Lambda_\alpha$ with those tied to
$\Lambda_\beta$, and involving the operator $^t\!\mathcal{X}_{\alpha,\beta}$.

\begin{Proposition}\label{proposition2}\qquad

\begin{enumerate}\itemsep=0pt
\item[$(i)$] The dual
Dunkl--Sonine integral transform $^t\!\mathcal{X}_{\alpha,\beta}$ maps  $L^1 (\mathbb{R}, |x|^{2\beta + 1}dx)$
continuously into $L^1 (\mathbb{R}, |x|^{2\alpha + 1}dx)$.

\item[$(ii)$] For every $f\in L^1 (\mathbb{R}, |x|^{2\beta + 1}dx)$ we have the identity \begin{gather}{\cal F}_\beta( f)
=\mathcal{F}_\alpha \circ {}^t\!\mathcal{X}_{\alpha,\beta}(f).\label{eq10}\end{gather}

\item[$(iii)$] Let $f, g\in L^1 (\mathbb{R}, |x|^{2\beta + 1}dx)$. Then
\begin{gather}^t\!\mathcal{X}_{\alpha,\beta}(f\ast_\beta g) = {}^t\!\mathcal{X}_{\alpha,\beta}f \ast_\alpha
{}^t\!\mathcal{X}_{\alpha,\beta}g.\label{eq11}
\end{gather}
\end{enumerate}
\end{Proposition}

\begin{proof}  Let $f \in L^1
(\mathbb{R}, |x|^{2\beta + 1}dx)$. By Fubini's theorem we have
\begin{gather*}
\int_{\mathbb{R}}{}^t\!\mathcal{X}_{\alpha,\beta}(|f|)(y)|y|^{2\alpha+1}dy = \int_{\mathbb{R}}\left(\int_{|x|\geq
|y|}\mathcal{K}_{\alpha,\beta}(x,y) |f(x)|
 |x|^{2\beta+1}\,dx\right)|y|^{2\alpha+1}dy\\
\phantom{\int_{\mathbb{R}}{}^t\!\mathcal{X}_{\alpha,\beta}(|f|)(y)|y|^{2\alpha+1}dy }{}  = \int_{\mathbb{R}}|f(x)|\left(\int_{-|x|}^{|x|}\mathcal{K}_{\alpha,\beta}(x,y)
|y|^{2\alpha+1}dy\right)|x|^{2\beta+1}\,dx.
\end{gather*}
But by \eqref{eq1},
\begin{gather}\int_{-|x|}^{|x|}\mathcal{K}_{\alpha,\beta}(x,y)
|y|^{2\alpha+1}dy=e_\beta(0)=1.\label{eq12}\end{gather}
Hence, $^t\!\mathcal{X}_{\alpha,\beta}f$ is almost
everywhere def\/ined on $\mathbb{R}$, belongs to $L^1 (\mathbb{R}, |x|^{2\alpha + 1}dx)$ and
$||^t\!\mathcal{X}_{\alpha,\beta}f||_{1,\alpha}\leq||f||_{1,\beta}$, which proves $(i)$. Identity \eqref{eq10} follows by
using \eqref{eq1}, \eqref{eq3}, \eqref{eq9}, and Fubini's theorem. Identity \eqref{eq11} follows by applying the Dunkl transform
$\mathcal{F}_\alpha$ to both its sides and by using \eqref{eq8}, \eqref{eq10} and Remark~\ref{remark1}.
\end{proof}

\begin{Remark} \label{remark3}
From \eqref{eq10} and
Remark~\ref{remark1}, we deduce that the transform
 $^t\!\mathcal{X}_{\alpha,\beta}$ maps $L^1
(\mathbb{R}$, $|x|^{2\beta + 1}dx)$ injectively into $L^1 (\mathbb{R}, |x|^{2\alpha + 1}dx)$.
\end{Remark}

From \cite{Soltani2008} we have the following result.

\begin{Theorem}\label{theorem3} The dual Dunkl--Sonine integral transform
$^t\!\mathcal{X}_{\alpha,\beta}$ is an automorphism of $\mathcal{S}(\mathbb{R})$ sa\-tisfying the intertwining
relation
\[^t\!\mathcal{X}_{\alpha,\beta} \Lambda_{\beta} f=\Lambda_{\alpha}{}^t\!\mathcal{X}_{\alpha,\beta} f, \qquad
f\in \mathcal{S}(\mathbb{R}).\]
Moreover $^t\!\mathcal{X}_{\alpha,\beta}$ admits the factorization
\[^t\!\mathcal{X}_{\alpha,\beta}f={}^tV_{\alpha}^{-1}\circ{} ^tV_{\beta} f\qquad \mbox{for  all} \ \  f\in
\mathcal{S}(\mathbb{R}),\]
where for $\gamma >-1/2$, $^tV_\gamma$ denotes the dual Dunkl intertwining operator
given by \[^tV_{\gamma}f(y)=\frac{\Gamma(\gamma+1)}{\sqrt{\pi}\,\Gamma(\gamma
+1/2)} \int_{|x|\geq|y|}\mbox{\rm sgn}(x)\,(x+y)\left(x^2-y^2\right)^{\gamma -1/2} f(x) \,dx .\] \end{Theorem}

\begin{Definition}[see \cite{Soltani2008}] \label{definition2} The Dunkl--Sonine integral transform $\mathcal{X}_{\alpha,\beta}$ is def\/ined for a
locally bounded function $f$ on $\mathbb{R}$ by
\begin{gather}
\mathcal{X}_{\alpha,\beta}f(x)=\left\{
\begin{array}{ll}
 \displaystyle{\int_{-|x|}^{|x|}}\mathcal{K}_{\alpha,\beta}(x,y)\,f(y)
\,|y|^{2\alpha+1}\,dy
 &\mbox{ if } x\neq 0,\\\\
 f(0)&\mbox{ if } x=0.
 \end{array}\right.\label{eq13}\end{gather}
  \end{Definition}

\begin{Remark}\label{remark4}
\qquad
\begin{enumerate}\itemsep=0pt
\item[$(i)$] Notice that by \eqref{eq12}, $||\mathcal{X}_{\alpha,\beta}f||_\infty\leq||f||_\infty$ if $f\in
L^\infty(\mathbb{R})$.

\item[$(ii)$] It follows from \eqref{eq1} that \begin{gather}e_\beta(\lambda x)=\mathcal{X}_{\alpha,\beta}(e_\alpha(\lambda
\cdot)(x)\label{eq14}\end{gather}for
 all $\lambda\in\mathbb{C}$ and $x\in \mathbb{R}$.
\end{enumerate}
\end{Remark}

\begin{Proposition} \label{proposition3} \qquad

\begin{enumerate}\itemsep=0pt
\item[$(i)$] For any $f\in L^\infty(\mathbb{R})$ and $g\in L^1 (\mathbb{R}, |x|^{2\beta + 1}dx)$ we
have the duality relation
\begin{gather}\int_{\mathbb{R}}\mathcal{X}_{\alpha,\beta}f(x)g(x)|x|^{2\beta + 1}dx=
\int_{\mathbb{R}}f(y)\,{}^t\!\mathcal{X}_{\alpha,\beta}g(y) |y|^{2\alpha+1}dy.\label{eq15}\end{gather}

\item[$(ii)$] Let $f\in L^1 (\mathbb{R}, |x|^{2\beta + 1}dx)$ and $g\in L^\infty(\mathbb{R})$. Then
\begin{gather}\mathcal{X}_{\alpha,\beta}\big(  {}^t\!\mathcal{X}_{\alpha,\beta}  f\ast_\alpha  g  \big)=f \ast_\beta
\mathcal{X}_{\alpha,\beta}g.\label{eq16}\end{gather}
\end{enumerate}
\end{Proposition}

 \begin{proof} Identity \eqref{eq15} follows by
using \eqref{eq9}, \eqref{eq13} and Fubini's theorem. Let us check \eqref{eq16}. Let $\psi\in \mathcal{S}(\mathbb{R})$. By using \eqref{eq11},
\eqref{eq15} and Fubini's theorem, we have
\begin{gather*}
\int_{\mathbb{R}}f \ast_\beta
\mathcal{X}_{\alpha,\beta}g(x) \psi(x)|x|^{2\beta+1}dx = \int_{\mathbb{R}}\mathcal{X}_{\alpha,\beta}g(x)
 \psi\!\ast_\beta \!f^-(x)\,|x|^{2\beta+1}dx\\
 \qquad{}= \int_{\mathbb{R}}g(y)\,{}^t\!\mathcal{X}_{\alpha,\beta}(\psi\!\ast_\beta \!f^-)(y) |y|^{2\alpha+1}dy
 = \int_{\mathbb{R}}g(y)\left( {}^t\!\mathcal{X}_{\alpha,\beta}\psi\ast_\alpha
 {}^t\!\mathcal{X}_{\alpha,\beta}f^-\right) (y) |y|^{2\alpha+1}dy,
 \end{gather*}
 where $f^-(x)=f(-x)$, $x\in\mathbb{R}$. But an easy computation shows that $^t\!\mathcal{X}_{\alpha,\beta}f^-
 =\left(^t\!\mathcal{X}_{\alpha,\beta}f\right)^-$. Hence,
\begin{gather*}
\int_{\mathbb{R}}f \ast_\beta \mathcal{X}_{\alpha,\beta}g(x) \psi(x)|x|^{2\beta+1}dx
  = \int_{\mathbb{R}}g(y)\, ^t\!\mathcal{X}_{\alpha,\beta}\psi\ast_\alpha
 \left({}^t\!\mathcal{X}_{\alpha,\beta}f\right)^- (y) |y|^{2\alpha+1}dy\\
 \qquad{}= \int_{\mathbb{R}}{}^t\!\mathcal{X}_{\alpha,\beta}f\ast_\alpha \!g(y)\,^t\!\mathcal{X}_{\alpha,\beta}\psi(y) |y|^{2\alpha+1}dy
  = \int_{\mathbb{R}}\mathcal{X}_{\alpha,\beta}\left(^t\!\mathcal{X}_{\alpha,\beta}f\ast_\alpha g\right) (x) \psi(x) |x|^{2\beta+1}dx.
\end{gather*} This clearly yields the result.\end{proof}

\section{Calder\'{o}n's formula for the Dunkl continuous\\ wavelet transform}\label{section4}
Throughout this section assume $\gamma
> - 1/2$.
\begin{Definition}\label{definition3}
 We say that a function  $g\in L^2(\mathbb{R},|x|^{2\gamma+1}dx)$ is a Dunkl wavelet of order
$\gamma$, if it satisf\/ies the admissibility
condition\begin{gather}0<C_g^\gamma:=\int_0^{\infty}|\mathcal{F}_\gamma
g(\lambda)|^2\frac{d\lambda}{\lambda}=\int_0^{\infty}|\mathcal{F}_\gamma
g(-\lambda)|^2\frac{d\lambda}{\lambda}<\infty.\label{eq17}\end{gather}
\end{Definition}
\begin{Remark} \label{remark5}\qquad
\begin{enumerate}\itemsep=0pt
\item[$(i)$] If
$g$ is real-valued we have $\mathcal{F}_\gamma g(-\lambda)=\overline{\mathcal{F}_\gamma g(\lambda)}$, so \eqref{eq17}
reduces to
\[0<C_g^\gamma:=\int_0^{\infty}|\mathcal{F}_\gamma g(\lambda)|^2\frac{d\lambda}{\lambda}<\infty.\]

\item[$(ii)$] If $0\neq g\in L^2(\mathbb{R},|x|^{2\gamma+1}dx)$ is real-valued and
satisf\/ies\[\exists\,\eta>0\qquad\mbox{such that}\quad \mathcal{F}_\gamma g(\lambda)-\mathcal{F}_\gamma
g(0)=\mathcal{O}(\lambda^\eta)\quad \mbox{as}\quad \lambda\rightarrow 0^+\]
then \eqref{eq17} is equivalent to
$\mathcal{F}_\gamma g(0)=0$.
\end{enumerate}
\end{Remark}

\begin{Note}\label{note2} For a function $g$ in $ L^2(\mathbb{R},|x|^{2\gamma+1}dx)$
and for $(a,b)\in (0,\infty)\times \mathbb{R}$ we
write\[g_{a,b}^\gamma(x):=\frac{1}{a^{2\gamma+2}} T_\gamma^{-b}g_a(x),\]
where $T_\gamma^{-b}$ are the generalized
translation operators given by \eqref{eq5}, and $g_a(x):=g(x/a)$, $x\in \mathbb{R}$. \end{Note}

\begin{Remark} \label{remark6} Let $g\in
L^2(\mathbb{R},|x|^{2\gamma+1}dx)$ and $a>0$. Then it is easily checked that $g_a\in
L^2(\mathbb{R}$, $|x|^{2\gamma+1}dx)$,
 $\left|\left|g_a\right|\right|_{2,\gamma}=a^{\gamma+1}\left|\left|g\right|\right|_{2,\gamma}$,
 and
 $\mathcal{F}_\gamma(g_a)(\lambda)=a^{2\gamma+2}\mathcal{F}_\gamma(g)(a\lambda)$.
 \end{Remark}

\begin{Definition}\label{definition4}
Let $g\in L^2(\mathbb{R},|x|^{2\gamma+1}dx)$ be a Dunkl wavelet of order $\gamma$. We
def\/ine for regular functions $f$ on $\mathbb{R}$, the Dunkl continuous wavelet transform by \[\Phi_g^\gamma
(f)(a,b):=\int_{\mathbb{R}}f(x) \overline{g_{a,b}^\gamma(x)} |x|^{2\gamma+1}dx\]which can also be written in the
form \[\Phi_g^\gamma (f)(a,b)=\frac{1}{a^{2\gamma+2}}\,f\ast _\gamma \widetilde{g}_a(b),\]where $\ast_\gamma$ is
the generalized convolution product given by~\eqref{eq6}, and $\widetilde{g}_a(x):=\overline{g(-x/a)}$, $x\in
\mathbb{R}$.\end{Definition}

The Dunkl continuous wavelet transform has been investigated in depth in~\cite{Jouini} in which precise def\/initions,
examples, and a more complete discussion of its properties can be found. We look here for a Calder\'{o}n formula
for this transform. We start with some technical lemmas.
\begin{Lemma}\label{lemma1}
For all $f,g\in
L^2(\mathbb{R},|x|^{2\gamma+1}dx)$ and all $\psi \in \mathcal{S}(\mathbb{R})$ we have the identity
\[
\int_{\mathbb{R}}f\ast_\gamma g(x) \mathcal{F}^{-1}_\gamma\psi(x)
  |x|^{2\gamma+1}dx=m_\gamma\int_{\mathbb{R}}\mathcal{F}_\gamma f(\lambda) \mathcal{F}_\gamma
  g(\lambda) \psi^-(\lambda) |\lambda|^{2\gamma+1}d\lambda,
  \]where $m_\gamma$ is given by~\eqref{eq4}.\end{Lemma}

  \begin{proof} Fix $g\in L^2(\mathbb{R},|x|^{2\gamma+1}dx)$ and $\psi\in \mathcal{S}(\mathbb{R})$.
   For $f\in L^2(\mathbb{R},|x|^{2\gamma+1}dx)$ put\[S_1(f):=\int_{\mathbb{R}}f\ast_\gamma g(x) \mathcal{F}^{-1}_\gamma\psi(x)
   |x|^{2\gamma+1}dx\] and
   \[
  S_2(f):=m_\gamma\int_{\mathbb{R}}\mathcal{F}_\gamma f(\lambda) \mathcal{F}_\gamma
  g(\lambda) \psi^-(\lambda) |\lambda|^{2\gamma+1}d\lambda.\]By \eqref{eq7}, \eqref{eq8} and Theorem~\ref{theorem1}, we see that $S_1(f)=S_2(f)$ for each $f\in L^1\cap L^2(\mathbb{R},|x|^{2\gamma+1}dx)$. Moreover, by using~\eqref{eq7},
H\"{o}lder's inequality and Theorem~\ref{theorem1} we have
  \[|S_1(f)| \leq ||f \ast_\gamma g||_{\infty}  ||\mathcal{F}^{-1}_\gamma \psi||_{1,\gamma}\leq 4
  ||f||_{2,\gamma} ||g||_{2,\gamma} ||\mathcal{F}^{-1}_\gamma
  \psi||_{1,\gamma}\]and
    \begin{gather*}
    |S_2(f)|  \leq  m_\gamma\,||\mathcal{F}_\gamma f  \mathcal{F}_\gamma
g||_{1,\gamma}
  ||\psi||_{\infty} \leq  m_\gamma\,||\mathcal{F}_\gamma f||_{2,\gamma} ||\mathcal{F}_\gamma g||_{2,\gamma}
  ||\psi||_{\infty} = ||f||_{2,\gamma} || g||_{2,\gamma}
  ||\psi||_{\infty},
    \end{gather*}
which shows that the linear functionals $S_1$ and $S_2$ are bounded on $L^2(\mathbb{R},|x|^{2\gamma+1}dx)$.
Therefore $S_1 \equiv S_2$, and the lemma is proved.\end{proof}

 \begin{Lemma}\label{lemma2} Let $f_1,f_2\in L^2(\mathbb{R},|x|^{2\gamma+1}dx)$. Then $f_1\ast_\gamma
 f_2\in
 L^2(\mathbb{R},|x|^{2\gamma+1}dx)$
  if and only if $\mathcal{F}_\gamma f_1 \mathcal{F}_\gamma f_2\in L^2(\mathbb{R},|x|^{2\gamma+1}dx)$ and we have
  \[\mathcal{F}_\gamma(f_1\ast_\gamma f_2 )=\mathcal{F}_\gamma f_1 \mathcal{F}_\gamma f_2\]in the
  $L^2$-case.\end{Lemma}

  \begin{proof} Suppose $f_1\ast_\gamma
 f_2\in
 L^2(\mathbb{R},|x|^{2\gamma+1}dx)$. By Lemma~\ref{lemma1} and Theorem~\ref{theorem1}, we have for any $\psi\in\mathcal{S}(\mathbb{R})$,
     \begin{gather*}
     m_\gamma\int_{\mathbb{R}}\mathcal{F}_\gamma f_1(\lambda)\mathcal{F}_\gamma
  f_2(\lambda) \psi(\lambda) |\lambda|^{2\gamma+1}d\lambda = \int_{\mathbb{R}}f_1\ast_\gamma f_2(x) \mathcal{F}^{-1}_\gamma\psi^-(x)
  |x|^{2\gamma+1}dx\\
  \qquad{} = \int_{\mathbb{R}}f_1\ast_\gamma f_2(x) \overline{\mathcal{F}^{-1}_\gamma\overline{\psi}(x)}
   |x|^{2\gamma+1}dx = m_\gamma\int_{\mathbb{R}}\mathcal{F}_\gamma (f_1\ast_\gamma f_2)(\lambda) \psi(\lambda)
  |\lambda|^{2\gamma+1}d\lambda,
    \end{gather*}
    which shows that $\mathcal{F}_\gamma f_1 \mathcal{F}_\gamma f_2=\mathcal{F}_\gamma (f_1\ast_\gamma f_2)$.
    Conversely, if $\mathcal{F}_\gamma f_1 \mathcal{F}_\gamma f_2 \in L^2(\mathbb{R},|x|^{2\gamma+1}dx)$, then by
Lemma~\ref{lemma1} and Theorem~\ref{theorem1}, we have for any $\psi\in\mathcal{S}(\mathbb{R})$,
\begin{gather*}
\int_{\mathbb{R}}f_1\ast_\gamma f_2(x) \mathcal{F}^{-1}_\gamma\psi(x)
   |x|^{2\gamma+1}dx = m_\gamma\int_{\mathbb{R}}\mathcal{F}_\gamma f_1(\lambda) \mathcal{F}_\gamma
  f_2(\lambda) \overline{\widetilde{\psi}(\lambda)} |\lambda|^{2\gamma+1}d\lambda\\
  \qquad{} = \int_{\mathbb{R}}
  \mathcal{F}_\gamma^{-1}(\mathcal{F}_\gamma f_1 \mathcal{F}_\gamma
  f_2)(x) \mathcal{F}_\gamma^{-1}\psi(x) |x|^{2\gamma+1}dx,
\end{gather*}
which shows, in view of Theorem~\ref{theorem2}, that $f_1\ast_\gamma f_2=\mathcal{F}_\gamma^{-1}(\mathcal{F}_\gamma
f_1 \mathcal{F}_\gamma
  f_2)$. This achieves the proof of Lemma~\ref{lemma2}.
 \end{proof}

 A combination of Lemma~\ref{lemma2} and Theorem~\ref{theorem1} gives us the following.
  \begin{Lemma}\label{lemma3} Let $f_1,f_2\in L^2(\mathbb{R},|x|^{2\gamma+1}dx)$. Then\[\int_{\mathbb{R}}|f_1\ast_\gamma f_2(x)|^2
  |x|^{2\gamma+1}dx=m_\gamma\int_{\mathbb{R}}|\mathcal{F}_\gamma f_1(\lambda)|^2|\mathcal{F}_\gamma
  f_2(\lambda)|^2 |\lambda|^{2\gamma+1}d\lambda,
  \]where both sides are finite or infinite.
  \end{Lemma}

   \begin{Lemma}\label{lemma4} Let $g\in L^2(\mathbb{R},|x|^{2\gamma+1}dx)$ be a Dunkl wavelet of order $\gamma$ such
   that
   $\mathcal{F}_\gamma g
    \in L^\infty(\mathbb{R})$. For
   $0<\varepsilon<\delta<\infty$ define
   \begin{gather}G_{\varepsilon,\delta}(x):=\frac{1}{C_g^\gamma}\int_\varepsilon^\delta g_a\ast_\gamma \widetilde{g}_a(x)
   \frac{da}{a^{4\gamma+5}}\label{eq18}\end{gather}
   and\begin{gather}K_{\varepsilon,\delta}(\lambda):=\frac{1}{C_g^\gamma}\int_\varepsilon^\delta |\mathcal{F}_\gamma g(a\lambda)|^2
   \frac{da}{a}.\label{eq19}\end{gather}Then \begin{gather}G_{\varepsilon,\delta}\in L^2(\mathbb{R},|x|^{2\gamma+1}dx),\qquad
    K_{\varepsilon,\delta}\in \left(L^1 \cap L^2\right)(\mathbb{R},|x|^{2\gamma+1}dx),\label{eq20}\end{gather}and \[\mathcal{F}_\gamma (G_{\varepsilon,\delta})
    =K_{\varepsilon,\delta}.\]\end{Lemma}

    \begin{proof} Using Schwarz inequality for the measure $\displaystyle{\frac{da}{a^{4\gamma+5}}}$ we obtain\[|G_{\varepsilon,\delta}(x)|^2
    \leq \frac{1}{\left(C_g^\gamma\right)^2}\left(\int_\varepsilon^\delta\frac{da}{a^{4\gamma+5}}\right)
    \int_\varepsilon^\delta
     |g_a\ast_\gamma
    \widetilde{g}_a(x)|^2
    \frac{da}{a^{4\gamma+5}},\]
    so
 \[\int_{\mathbb{R}}|G_{\varepsilon,\delta}(x)|^2|x|^{2\gamma+1}dx\leq
    \frac{1}{\left(C_g^\gamma\right)^2}\left(\int_\varepsilon^\delta\frac{da}{a^{4\gamma+5}}\right)\int_\varepsilon^\delta
    \int_{\mathbb{R}}|g_a\ast_\gamma
    \widetilde{g}_a(x)|^2|x|^{2\gamma+1}dx
    \frac{da}{a^{4\gamma+5}}.
    \]
    By Theorem~\ref{theorem1}, Lemma~\ref{lemma3}, and Remark~\ref{remark6}, we have
  \begin{gather*}
    \int_{\mathbb{R}}|g_a\ast_\gamma
    \widetilde{g}_a(x)|^2|x|^{2\gamma+1}dx = m_\gamma\int_{\mathbb{R}}|\mathcal{F}_\gamma (g_a)(\lambda)|^
    4|\lambda|^{2\gamma+1}d\lambda\\
 \qquad{}  \leq m_\gamma\left|\left|\mathcal{F}_\gamma
  (g_a)\right|\right|_\infty^2\int_{\mathbb{R}}|\mathcal{F}_\gamma (g_a)(\lambda)|^2|\lambda|^{2\gamma+1}d\lambda \\
  \qquad{}= \left|\left|\mathcal{F}_\gamma
  (g_a)\right|\right|_\infty^2\left|\left|g_a\right|\right|_{2,\gamma}^2 = a^{6\gamma+6}\left|\left|\mathcal{F}_\gamma
  g\right|\right|_\infty^2\left|\left|g\right|\right|_{2,\gamma}^2.
  \end{gather*}
Hence \[\int_{\mathbb{R}}|G_{\varepsilon,\delta}(x)|^2|x|^{2\gamma+1}dx\leq
    \frac{\left|\left|\mathcal{F}_\gamma
  g\right|\right|_\infty^2\left|\left|g\right|\right|_{2,\gamma}^2}{\left(C_g^\gamma\right)^2}
  \left(\int_\varepsilon^\delta a^{2\gamma+1}da\right) \left(\int_\varepsilon^\delta\frac{da}{a^{4\gamma+5}}
  \right)<\infty.\]
The second assertion in \eqref{eq20} is easily checked. Let us calculate $\mathcal{F}_\gamma (G_{\varepsilon,\delta})$.
Fix $x\in \mathbb{R}$. From Theorem~\ref{theorem1} and Lemma~\ref{lemma2} we get\[g_a\ast_\gamma
\widetilde{g}_a(x)=m_\gamma\int_{\mathbb{R}}|\mathcal{F}_\gamma (g_a)(\lambda)|^2e_{\gamma}(i\lambda x)
|\lambda|^{2\gamma+1}d\lambda,\]so
\[G_{\varepsilon,\delta}(x)=\frac{m_\gamma}{C_g^\gamma}\int_\varepsilon^\delta\left(\int_{\mathbb{R}}|\mathcal{F}_\gamma
(g_a)(\lambda)|^2e_{\gamma}(i\lambda x) |\lambda|^{2\gamma+1}d\lambda\right)\frac{da}{a^{4\gamma+5}}.\] As
$|e_{\gamma}(iz)|\leq 1$ for all $z\in \mathbb{R}$ (see~\cite{Rosler}), we deduce by Theorem~\ref{theorem1} that
\begin{gather*}
m_\gamma\int_\varepsilon^\delta\int_{\mathbb{R}}|\mathcal{F}_\gamma (g_a)(\lambda)|^2|e_{\gamma}(i\lambda x)|
|\lambda|^{2\gamma+1}d\lambda\frac{da}{a^{4\gamma+5}}\\
\qquad{} \leq
\int_\varepsilon^\delta||g_a||_{2,\gamma}^2\frac{da}{a^{4\gamma+5}} = ||g||_{2,\gamma}^2\int_\varepsilon^\delta
\frac{da}{a^{2\gamma+3}}<\infty.
  \end{gather*}
Hence, applying Fubini's theorem, we f\/ind that
  \begin{gather*}
G_{\varepsilon,\delta}(x) = m_\gamma\int_{\mathbb{R}}\left(\frac{1}{C_g^\gamma}\int_\varepsilon^\delta|\mathcal{F}_\gamma
g(a\lambda)|^2\frac{da}{a}\right)e_{\gamma}(i\lambda x)
|\lambda|^{2\gamma+1}d\lambda\\
\phantom{G_{\varepsilon,\delta}(x)}{} = m_\gamma\int_{\mathbb{R}}K_{\varepsilon,\delta}(\lambda)e_{\gamma}(i\lambda x)
|\lambda|^{2\gamma+1}d\lambda
  \end{gather*}
which completes the proof.\end{proof}

We can now state the main result of this section.

\begin{Theorem}[Calder\'{o}n's formula]\label{theorem4}  Let
$g\in L^2(\mathbb{R},|x|^{2\gamma+1}dx)$ be a Dunkl wavelet of order $\gamma$ such that $\mathcal{F}_\gamma g \in
L^{\infty}(\mathbb{R})$. Then for $f \in L^2(\mathbb{R},|x|^{2\gamma+1}dx)$
  and $0<\varepsilon<\delta<\infty$, the function
  \[f^{\varepsilon,\delta}(x):=\frac{1}{C_g^\gamma}\int_\varepsilon^\delta\int_{\mathbb{R}} \Phi_g^\gamma
(f)(a,b)g_{a,b}(x)|b|^{2\gamma+1}db\frac{da}{a}\]belongs to $L^2(\mathbb{R},|x|^{2\gamma+1}dx)$ and
satisfies\begin{gather}\lim_{\varepsilon\rightarrow 0,\,\delta\rightarrow \infty}
\big\|f^{\varepsilon,\delta}-f\big\|_{2,\gamma}=0.\label{eq21}\end{gather}
\end{Theorem}

\begin{proof} It is easily seen that\[f^{\varepsilon,\delta}=f\ast_\gamma
G_{\varepsilon,\delta},
\]
where $G_{\varepsilon,\delta}$ is given by~\eqref{eq18}. It follows by Lemmas~\ref{lemma2} and~\ref{lemma4} that
$f^{\varepsilon,\delta}\in L^2(\mathbb{R},|x|^{2\gamma+1}dx)$ and
$\mathcal{F}_\gamma(f^{\varepsilon,\delta})=\mathcal{F}_\gamma(f)\,K_{\varepsilon,\delta}$,
 where $K_{\varepsilon,\delta}$ is as in~\eqref{eq19}. From this and
 Theorem~\ref{theorem1} we obtain
  \begin{gather*}
 \big\|f^{\varepsilon,\delta}-f\big\|_{2,\gamma}^2 = m_\gamma\int_{\mathbb{R}}|\mathcal{F}_\gamma(f^{\varepsilon,\delta}
 -f)(\lambda)|^2
 |\lambda|^{2\gamma+1}d\lambda\\
 \phantom{\big\|f^{\varepsilon,\delta}-f\big\|_{2,\gamma}^2}{} = m_\gamma\int_{\mathbb{R}}|\mathcal{F}_\gamma f(\lambda)|^2(1-
 K_{\varepsilon,\delta}(\lambda))^2
 |\lambda|^{2\gamma+1}d\lambda.
   \end{gather*}
But by~\eqref{eq17} we have\[\lim_{\varepsilon\rightarrow 0,\,\delta\rightarrow
\infty}K_{\varepsilon,\delta}(\lambda)=1,\qquad \mbox{for almost all}\ \ \lambda\in \mathbb{R}.\]
So \eqref{eq21} follows
from the dominated convergence theorem. \end{proof}

Another pointwise inversion formula for the Dunkl wavelet transform, proved in~\cite{Jouini},  is as
follows.

\begin{Theorem}\label{theorem5} Let $g\in L^2(\mathbb{R},|x|^{2\gamma+1}dx)$ be a Dunkl wavelet of order
$\gamma$. If both $f$ and
 $\mathcal{F}_\gamma f$ are in $L^1(\mathbb{R},|x|^{2\gamma+1}dx)$ then we have
  \[f(x)=\frac{1}{C_g^\gamma}\int_0^{\infty}\left(\int_{\mathbb{R}} \Phi_g^\gamma
(f)(a,b)g_{a,b}^\gamma(x)|b|^{2\gamma+1}db\right)\frac{da}{a},\quad \mbox{a.e.,} \]where, for each $x\in
\mathbb{R}$, both the inner integral and the outer integral are absolutely convergent, but possibly not the double
integral.
\end{Theorem}

\section[Inversion of the dual
Dunkl-Sonine transform  using Dunkl wavelets]{Inversion of the dual
Dunkl--Sonine transform\\ using Dunkl wavelets}\label{section5}

From now on assume $\beta>\alpha>-1/2$. In order to invert the dual Dunkl--Sonine transform, we need the following
two technical lemmas.

\begin{Lemma}\label{lemma5} Let $0\neq g\in L^1\cap L^2(\mathbb{R},|x|^{2\alpha+1}dx)$ such
that
 $\mathcal{F}_\alpha g \in L^1(\mathbb{R},|x|^{2\alpha+1}dx)$ and satisfying
 \begin{gather}\exists \, \eta >\beta-2\alpha-1 \qquad \mbox{such that}\quad \mathcal{F}_\alpha g(\lambda)=\mathcal{O}\left(|\lambda|^\eta\right)
 \quad \mbox{as}\quad\lambda\rightarrow 0.\label{eq22}\end{gather}
 Then $\mathcal{X}_{\alpha,\beta}g \in L^2(\mathbb{R},|x|^{2\beta+1}dx)$ and
 \[\mathcal{F}_\beta(\mathcal{X}_{\alpha,\beta}g)(\lambda)
 =\frac{m_\alpha}{m_\beta}\,\frac{\mathcal{F}_\alpha
g(\lambda)}
 {|\lambda|^{2(\beta-\alpha)}}.\]
 \end{Lemma}
 \begin{proof} By Theorem~\ref{theorem1} we have\[g(x)=m_\alpha\int_{\mathbb{R}}\mathcal{F}_\alpha g(\lambda)e_\alpha(i\lambda x)
 |\lambda|^{2\alpha+1}d\lambda,\qquad \mbox{a.e.}\]
 So using~\eqref{eq14}, we f\/ind
that
\begin{gather}\mathcal{X}_{\alpha,\beta}g(x)=m_\beta\int_{\mathbb{R}}h_{\alpha,\beta}(\lambda)e_\beta(i\lambda
x)
 |\lambda|^{2\beta+1}d\lambda,\qquad \mbox{a.e.}\label{eq23}\end{gather}
 with\[h_{\alpha,\beta}(\lambda):=\frac{m_\alpha}{m_\beta}\,\frac{\mathcal{F}_\alpha g(\lambda)}
 {|\lambda|^{2(\beta-\alpha)}}.\]
 Clearly, $h_{\alpha,\beta}\in L^1(\mathbb{R},|x|^{2\beta+1}dx)$. So it suf\/f\/ices, in view of~\eqref{eq23}  and Theorem~\ref{theorem1},
 to prove that
 $h_{\alpha,\beta}$ belongs to $L^2(\mathbb{R},|x|^{2\beta+1}dx)$. We have
\begin{gather*}
\int_{\mathbb{R}}|h_{\alpha,\beta}(\lambda)|^2|\lambda|^{2\beta+1}d\lambda = \left(\frac{m_\alpha}{m_\beta}\right)^2\int_{\mathbb{R}}|\lambda|^{4\alpha-2\beta+1}
 |\mathcal{F}_\alpha g(\lambda)|^2d\lambda\\
 \qquad{} = \left(\frac{m_\alpha}{m_\beta}\right)^2\left(\int_{|\lambda|\leq
 1}+\int_{|\lambda|\geq 1}\right)|\lambda|^{4\alpha-2\beta+1}
 |\mathcal{F}_\alpha g(\lambda)|^2d\lambda := I_1+I_2.
\end{gather*}

By~\eqref{eq22} there is a positive constant $k$ such that\[I_1\leq k \int_{|\lambda|\leq
 1}|\lambda|^{2\eta +4\alpha-2\beta+1}d\lambda=\frac{k}{\eta+2\alpha-\beta+1}<\infty.\]

 From Theorem~\ref{theorem1}, it follows that
 \begin{gather*}
I_2 = \left(\frac{m_\alpha}{m_\beta}\right)^2\int_{|\lambda|\geq 1}|\lambda|^{2(\alpha-\beta)}
 |\mathcal{F}_\alpha
 g(\lambda)|^2|\lambda|^{2\alpha+1}d\lambda\\
 \phantom{I_2}{} \leq \left(\frac{m_\alpha}{m_\beta}\right)^2\int_{|\lambda|\geq
 1}
 |\mathcal{F}_\alpha
 g(\lambda)|^2|\lambda|^{2\alpha+1}d\lambda \leq \left(\frac{m_\alpha}{m_\beta}\right)^2\left|\left|\mathcal{F}_\alpha
 g\right|\right|_{2,\alpha}^2=\frac{m_\alpha}{(m_\beta)^2}\,||g||_{2,\alpha}^2<\infty
 \end{gather*}
 which ends the proof.\end{proof}
 \begin{Lemma}\label{lemma6} Let $0\neq g\in\!
L^1\cap L^2(\mathbb{R},|x|^{2\alpha+1}dx)$ be real-valued such that $\mathcal{F}_\alpha g \in\!
L^1(\mathbb{R},|x|^{2\alpha+1}dx)$ and satisfying
 \begin{gather}\exists \, \eta >2(\beta-\alpha) \qquad \mbox{such that}\quad \mathcal{F}_\alpha g(\lambda)=\mathcal{O}(\lambda^\eta)
 \quad \mbox{as}\quad\lambda\rightarrow 0^+.\label{eq24}
 \end{gather}
 Then $\mathcal{X}_{\alpha,\beta}g \in L^2(\mathbb{R},|x|^{2\beta+1}dx)$ is a Dunkl wavelet of order $\beta$
 and
 $\mathcal{F}_\beta(\mathcal{X}_{\alpha,\beta}g)\in
 L^\infty(\mathbb{R})$.
\end{Lemma}

\begin{proof} By combining~\eqref{eq24} and Lemma~\ref{lemma5} we see that $\mathcal{X}_{\alpha,\beta}g \in
L^2(\mathbb{R},|x|^{2\beta+1}dx)$,
 $\mathcal{F}_\beta(\mathcal{X}_{\alpha,\beta}g)$ is
 bounded and \[\mathcal{F}_\beta(\mathcal{X}_{\alpha,\beta}g)(\lambda)=\mathcal{O}\big(\lambda^{\eta-2(\beta-\alpha)}\big)
 \qquad \mbox{as}\quad\lambda\rightarrow 0^+.\]Thus, in view of Remark~\ref{remark5}, $\mathcal{X}_{\alpha,\beta}g$
  satisf\/ies the admissibility condition~\eqref{eq17} for $\gamma=\beta$.
  \end{proof}

  \begin{Remark} \label{remark7} In view of Remark~\ref{remark5}, each function
  satisfying the conditions of Lemma~\ref{lemma5} is a~Dunkl wavelet of order
  $\alpha$.\end{Remark}

  \begin{Lemma}\label{lemma7} Let $g$ be as in Lemma~{\rm \ref{lemma6}}. Then for all $f\in L^1(\mathbb{R},|x|^{2\beta+1}dx)$ we have
  \[\Phi_{\mathcal{X}_{\alpha,\beta}g}^\beta (f)(a,b)=\frac{1}{a^{2(\beta-\alpha)}}\,
  \mathcal{X}_{\alpha,\beta}\big[\Phi_g^\alpha
  \left(^t\!\mathcal{X}_{\alpha,\beta}f\right)(a,\cdot)\big](b).\]\end{Lemma}

  \begin{proof} By Def\/inition~\ref{definition4} we have\[\Phi_{\mathcal{X}_{\alpha,\beta}g}^\beta
(f)(a,b)=\frac{1}{a^{2\beta+2}}\,f\ast _\beta \widetilde{\left(\mathcal{X}_{\alpha,\beta}g\right)_a}(b).\]But
$\widetilde{\left(\mathcal{X}_{\alpha,\beta}g\right)_a}=
  \mathcal{X}_{\alpha,\beta}\left(\widetilde{g}_a\right)$ by
  virtue of~\eqref{eq2} and~\eqref{eq13}. So using~\eqref{eq16} we f\/ind
  that
   \begin{gather*}
\Phi_{\mathcal{X}_{\alpha,\beta}g}^\beta (f)(a,b) = \frac{1}{a^{2\beta+2}} f\ast _\beta
\left[\mathcal{X}_{\alpha,\beta}\left(\widetilde{g}_a\right)\right](b)\\
\phantom{\Phi_{\mathcal{X}_{\alpha,\beta}g}^\beta (f)(a,b)}{} = \frac{1}{a^{2\beta+2}} \mathcal{X}_{\alpha,\beta}
\left[{}^t\!\mathcal{X}_{\alpha,\beta}f *_\alpha
\widetilde{g}_a\right](b) = \frac{1}{a^{2(\beta-\alpha)}} \mathcal{X}_{\alpha,\beta}\big[\Phi_g^\alpha
\left(^t\!\mathcal{X}_{\alpha,\beta}f\right)(a,\cdot)\big](b),
 \end{gather*}
 which gives the desired result.
\end{proof}

 Combining Theorems~\ref{theorem4},~\ref{theorem5} with Lemmas~\ref{lemma6},~\ref{lemma7} we
get\begin{Theorem}\label{theorem6} Let $g$ be as in Lemma~{\rm \ref{lemma6}}. Then we have the following inversion formulas for the
dual Dunkl--Sonine transform:

$(i)$ If both $f$ and $\mathcal{F}_\beta f$ are in $L^1(\mathbb{R},|x|^{2\beta+1}dx)$ then for almost all $x\in
\mathbb{R}$ we have \[ f(x)=\frac{1}{C_{\mathcal{X}_{\alpha,\beta}g}^\beta}\int_0^{\infty}\left(\int_{\mathbb{R}}
\mathcal{X}_{\alpha,\beta}\left[\Phi_g^\alpha
\left(\,^t\!\mathcal{X}_{\alpha,\beta}f\right)\!(a,\cdot)\right] (b) \big(\mathcal{X}_{\alpha,\beta}g\big)_{a,b}^\beta(x)
 |b|^{2\beta+1}db\right)\frac{da}{a^{2(\beta-\alpha)+1}}.\]

$(ii)$ For $f \in  L^1 \cap L^2(\mathbb{R},|x|^{2\beta+1}dx)$
  and $0<\varepsilon<\delta<\infty$, the function
  \[f^{\varepsilon,\delta}(x):=\frac{1}{C_{\mathcal{X}_{\alpha,\beta}g}^\beta}\int_\varepsilon^\delta \int_{\mathbb{R}}
\mathcal{X}_{\alpha,\beta}\left[\Phi_g^\alpha
\left(\,^t\!\mathcal{X}_{\alpha,\beta}f\right) (a,\cdot)\right] (b) \big(\mathcal{X}_{\alpha,\beta}g\big)_{a,b}^\beta(x)
 |b|^{2\beta+1}db\,\frac{da}{a^{2(\beta-\alpha)+1}}\] satisfies
 \[
 \lim_{\varepsilon\rightarrow
0,\,\delta\rightarrow \infty}
\big\|f^{\varepsilon,\delta}-f\big\|_{2,\beta}=0.
\]
\end{Theorem}
\subsection*{Acknowledgements}

The author is grateful to the referees and editors for careful reading and useful comments.

\pdfbookmark[1]{References}{ref}
\LastPageEnding

\end{document}